\documentclass[a4paper,12pt]{amsart}
\usepackage{amsmath,amsfonts,amssymb}
\usepackage{mathabx}
\usepackage{graphicx}
\usepackage{url}
\usepackage{hyperref}
\usepackage{comment}
\usepackage{makecell}
\usepackage{boldline}
\setcellgapes{3pt}
\usepackage{hhline,colortbl}
\usepackage{xcolor}

\calclayout
\theoremstyle{plain}
\newtheorem{theorem}{Theorem}[section]
\newtheorem{lemma}[theorem]{Lemma}

\theoremstyle{definition}

\theoremstyle{remark}

\newcommand{\Z}{\mathbb{Z}}

\definecolor{blueviolet}{rgb}{0.5,0,1}
\begin{document}

\title{Golden ratio on nonorientable surfaces}

\author{Ji-Young Ham and Joongul Lee}
\address{Da Vinci College of General Education, Chung-Ang University, 
84 HeukSeok-Ro, DongJak-Gu, Seoul,
06974,\\
  Korea} 
\email{jiyoungham1@gmail.com.}

\address{Department of Mathematics Education, Hongik University, 
94 Wausan-ro, Mapo-gu, Seoul,
04066\\
   Korea} 
\email{jglee@hongik.ac.kr}



\maketitle 

\markboth{Ji-Young Ham and Joongul Lee} 
{Golden ratio on nonorientable surfaces}

\date{\today}

\begin{abstract}
On each nonorientable surface of odd genus $g \geq 5$, we give a mapping class whose dilatation on an invariant subsurface is the golden ratio.
\end{abstract}

\maketitle
\section{Introduction} 

In 1970's, Thurston classified the mapping class group of a surface into periodic, pseudo-Anosov, and reducible~\cite{Thurston88}. The same content from a somewhat different point of view can be found in~\cite{CassonBleiler88}. In this paper, we adopt the Penner's approach~\cite{Penner88,Penner91}. Penner made use of bigon tracks, a slight generalization of train track. Nice example of bigon tracks can be found in~\cite{LiechtiStrenner18-1}. Our paper is based on 
Liechti-Strenner's pseudo-Anosov diffeomorphisms on nonorientable surfaces~\cite{LiechtiStrenner18}.

Let $\Sigma_g$ be a surface of finite type.
A homeomorphism $h$ of $\Sigma_g$ is called \emph{pseudo-Anosov} if there is a pair of transversely measured foliations 
$\mathcal{F}^u$ and $\mathcal{F}^s$ in $\Sigma_g$ and a real number $\lambda >1$ such that $h(\mathcal{F}^u)=\lambda \mathcal{F}^u$ and  $h(\mathcal{F}^s)=1/\lambda \mathcal{F}^s$~\cite{Thurston88,CassonBleiler88}.
The number $\lambda$ is called the \emph{dilatation} of $h$ and the logarithm of $\lambda$ is called the \emph{topological entropy}. 

A mapping class $\Phi$ on $\Sigma_g$ is reducible if there is a family of essential disjoint simple closed nonboundary- or puncture-parallel curves $C$ and a representative $h$ of $\Phi$ with $h(C) = C$. In case $S$ is a (not necessarily connected) $h$-component of $\Sigma_g-C$, $h |_{S}$ is either periodic or pseudo-Anosov~\cite[Theorem 4]{Thurston88}.
The main purpose of the paper is to  give a reducible mapping class on each nonorientable surface of odd genus $g \geq 5$ such that the dilatation on an invariant subsurface is the golden ratio. Our main theorem is stated in Theorem~\ref{thm:main}

\begin{theorem} \label{thm:main}
Let $k \geq 3$ be an odd natural number.
On the Liechti-Strenner surface $\Sigma_{2k,k}$ of genus $k+2$, the mapping class $\Phi_k=r \circ T_{c_1} \circ r^{k-1}$ has an invariant subsurface  
such that the dilatation on the invariant subsurface is the golden ratio.
\end{theorem}

Note that the invariant subsurface is a punctured torus and the restriction of the map to the subsurface is an orientation-reversing mapping class, namely the one given by the matrix $\begin{bmatrix} 0 & 1 \\ 1 & 1 \end{bmatrix}$ on the torus (This is actually the orientation-reversing Anosov mapping class with the smallest dilatation on the torus).
Note also that one could of course realize the golden ratio as the dilatation of a reducible mapping class of any nonorientable surface of genus $g$ at least 3 by just gluing the appropriate number of copies of real projective planes along the boundary of the punctured torus and extending by the identity. But this approach doesn't seem to reveal some kind of nice symmetries and full dynamics.

The proof of Theorem~\ref{thm:main} is given in Section~\ref{sec:proof}. The proof for $\Sigma_{6,3}$ and the proof for $\Sigma_{2k,k}$ are the essentially the same. We give the proof of $\Sigma_{6,3}$ as a warming up.
\section{Liechti-Strenner construction of nonorientable surfaces} \label{section:LS}

We will briefly introduce the Liechti-Strenner's method of constructing the nonorientable surface $\Sigma_{2k,k}$ of genus $g=k+2$.

\subsection{The graph $G_{2k,k}$}
\begin{figure}
\begin{center}
\includegraphics[scale=1]{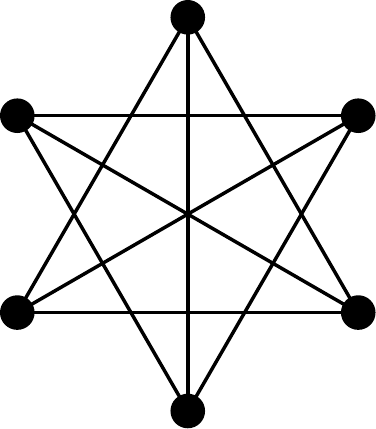}
\caption{The graph $G_{6,3}$.}
\label{fig:graph}
\end{center}
\end{figure}

Let $k \geq 3$ be an odd natural number. Let $G_{2k,k}$ be the graph whose vertices are the vertices of a regular $2k$-gon and every
vertex $v$ is connected to the $k$ vertices that are the farthest away from $v$ in the cyclic order of the vertices. Figure~\ref{fig:graph} shows the graph 
$G_{6,3}$.

\subsection{The surface $\Sigma_{2k,k}$}
\begin{figure}
\begin{center}
\includegraphics[scale=1]{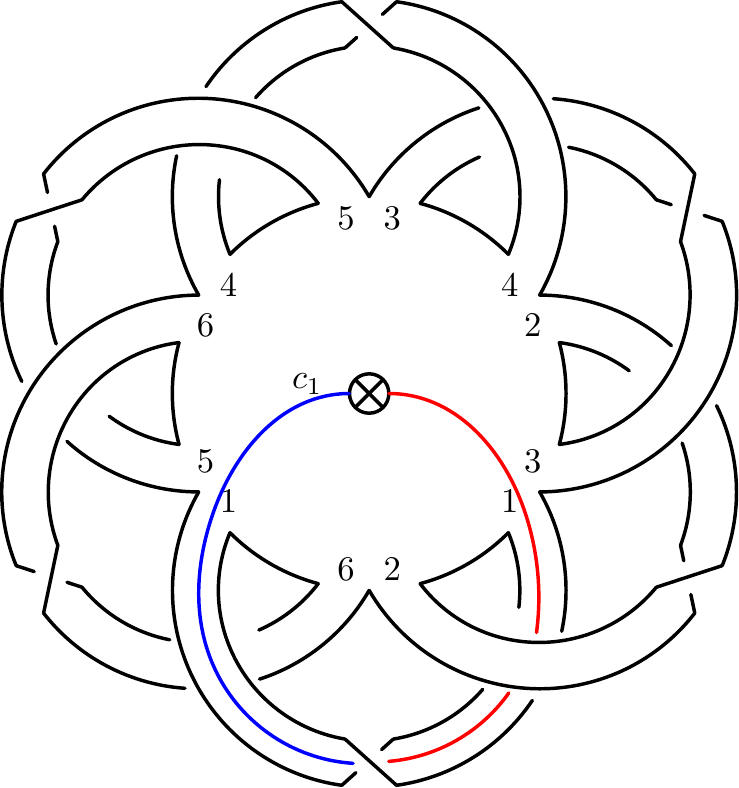}
\caption{The surface $\Sigma_{6,3}$ and the curve $c_1$.}
\label{fig:surface}
\end{center}
\end{figure}

For each $G_{2k,k}$, Liechti and Strenner constructed an nonorientable surface $\Sigma_{2k,k}$ that contains a collection of curves with
intersection graph $G_{2k,k}$~\cite[Subsection 2.2]{LiechtiStrenner18}. 
To construct $\Sigma_{2k,k}$, start with a disk with one crosscap. Next, we consider $4k$ disjoint intervals on the boundary of the disk and label
the intervals with integers from $1$ to $2k$ so that each label is used exactly twice.
In the cyclic order, the labels are $1,s,2,s+1,...,2k,s+2k$ where $s=\frac{3k+3}{2}$ and
all labels are understood modulo $2k$.
For each label, the corresponding two intervals are connected by a twisted strip. Figure~\ref{fig:surface} shows the surface $\Sigma_{6,3}$.

\begin{lemma}~\cite[Proposition 2.3]{LiechtiStrenner18}
The surface $\Sigma_{2k,k}$ is homeomorphic to the nonorientable surface of genus $k+2$ with
$k$ boundary components.
\end{lemma}

\subsection{The curves}
\begin{figure}
\begin{center}
\includegraphics[scale=1]{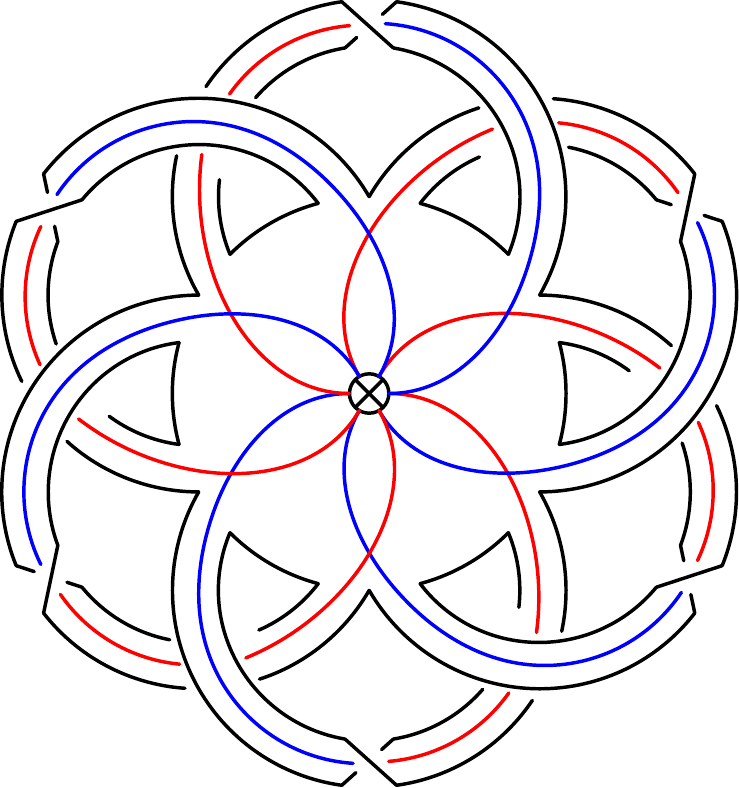}
\caption{A collection of filling inconsistently marked curves.}
\label{fig:curve}
\end{center}
\end{figure}

Liechti and Strenner constructed a two-sided curve $c_i$ for each label $i=1,\ldots,2k$ as follows. Each curve consists of two parts.
One part of each curve is the core of the strip corresponding to the label. The other part is an arc inside the disk that passes through the crosscap and connects the corresponding two intervals. The curve $c_1$ is shown in Figure~\ref{fig:surface}.
We choose markings for the $c_i$ which are invariant under the rotational symmetry See Figure~\ref{fig:curve}.

Note that every pair of curves intersects either once or not at all. The curves $c_i$ and $c_j$ are disjoint if and only if the two $i$ labels and the two $j$ labels link in the cyclic order. In other words, if the two $i$ labels separate the two $j$ labels.

\section {The mapping classes}

Denote by $r$ the rotation of $\Sigma_{2k,k}$ by one click in the clockwise direction. Denote by $T_{c_1}$ the right-handed Dehn twist about the curve $c_1$.
Define the mapping class
$$\Phi_k=r \circ T_{c_1} \circ r^{k-1}.$$

According to Penner~\cite[Theorem 3.1, Theorem 4.1]{Penner88} and Liechti and Strenner~\cite[Proposition 2.6]{LiechtiStrenner18}, we can construct a bigon track on 
$\Sigma_{2k,k}$ with filling curves $c_i$, $i=1,\ldots, 2k$. Each $c_i$ defines a characteristic measure $\mu_i$ on this bigon track, defined by assigning $1$ to the branches traversed by $c_i$ and zero to the rest. Let $H$ denote the cone generated by the measures $\mu_i$ in the cone of measures on the bigon track. 
Let $\mathcal{C}_k$ be $\mathcal{C}_k=\{c_i\,|i=1,\ldots,2k\}$ and $\mathcal{S}_k$ be the semigroup with presentation
$$\mathcal{S}_k (\mathcal{C}_k)=\langle c_i \in \mathcal{C}_k \, : \, c_i \leftrightarrow c_j \text{ if } c_i \cap c_j = \emptyset \rangle.$$

\vspace{2ex}

\begin{theorem}~\cite[Theorem 3.4]{Penner88}
The action of $\mathcal{S}_k (\mathcal{C}_k)$ on $H$ admits a faithful representation as a semigroup of invertible (over $\Z$) positive matrices.
\end{theorem}

\section{proof of Theorem~\ref{thm:main}} \label{sec:proof}
Note that each $c_i$ corresponds to $T_{c_1}$ and the action of $c_1$ on $H$ in the basis $\mu_i$ is given by
$$I_{2k}+A$$
where $I_{2k}$ is the $2k \times 2k$-identity matrix, and $A=\begin{bmatrix} a_{ij} \end{bmatrix}$ with $a_{ij}=0$ if $i \neq 1$, and $a_{1j}=\text{card}(c_1 \cap c_j)$.
The rotation $r$ acts by a permutation matrix.

\subsection{proof of the Theorem~\ref{thm:main} on $\Sigma_{6,3}$}
On $\Sigma_{6,3}$, the action of $c_1$ on $H$ in the basis $\{\mu_i\}$ is given by

$$\begin{bmatrix}
1 & 0 & 1 & 1 & 1 & 0 \\
0 & 1 & 0 & 0 & 0 & 0 \\
0 & 0 & 1 & 0 & 0 & 0 \\
0 & 0 & 0 & 1 & 0 & 0 \\
0 & 0 & 0 & 0 & 1 & 0 \\
0 & 0 & 0 & 0 & 0 & 1 \\
\end{bmatrix}.$$

The action of $r$ on $H$ in the basis $\mu_i$ is given by

$$\begin{bmatrix}
0 & 1 & 0 & 0 & 0 & 0 \\
0 & 0 & 1 & 0 & 0 & 0 \\
0 & 0 & 0 & 1 & 0 & 0 \\
0 & 0 & 0 & 0 & 1 & 0 \\
0 & 0 & 0 & 0 & 0 & 1 \\
1 & 0 & 0 & 0 & 0 & 0 \\
\end{bmatrix}.$$

Hence $\Phi_3$ is represented by

$$\begin{bmatrix}
0 & 0 & 0 & 1 & 0 & 0 \\
0 & 0 & 0 & 0 & 1 & 0 \\
0 & 0 & 0 & 0 & 0 & 1 \\
1 & 0 & 0 & 0 & 0 & 0 \\
0 & 1 & 0 & 0 & 0 & 0 \\
1 & 0 & 1 & 0 & 1 & 1 \\
\end{bmatrix}.$$

By changing the order of basis by switching $\mu_3$ and $\mu_5$, $\Phi_3$ can be represented by the following reducible matrix $M_3$,

$$\begin{bmatrix}
0 & 0 & 0 & 1 & 0 & 0 \\
0 & 0 & 1 & 0 & 0 & 0 \\
0 & 1 & 0 & 0 & 0 & 0 \\
1 & 0 & 0 & 0 & 0 & 0 \\
0 & 0 & 0 & 0 & 0 & 1 \\
1 & 0 & 1 & 0 & 1 & 1 \\
\end{bmatrix}.$$

$M_3$ can be partitioned into four submatrices. The partitioned matrix can be written as
$$\begin{bmatrix}
M_{11} & \mathbf{0} \\
M_{21}& M_{22} \\
\end{bmatrix}$$

where $M_{11}$ is the $4 \times 4$ matrix, $\mathbf{0}$ is the $4 \times 2$ matrix, $M_{21}$ is the $2 \times 4$ matrix, and
$M_{22}$ is the $2 \times 2$ matrix.

The characteristic polynomial of $M_3$ is 
\begin{align*}
\text{det}(xI-M_{3} )&=
\text{det}
\begin{bmatrix}
xI-M_{11} & \mathbf{0} \\
\mathbf{0} & I \\
\end{bmatrix}
\text{det}
\begin{bmatrix}
I & \mathbf{0}\\
-M_{21} & I \\
\end{bmatrix}
\text{det}
\begin{bmatrix}
I & \mathbf{0} \\
\mathbf{0} & xI-M_{22}\\
\end{bmatrix} \\
&=\text{det}(xI-M_{11} )\text{det}(xI-M_{22} ) \\
&=(x+1)^2(x-1)^2(x^2-x-1).
\end{align*}

\begin{figure}
\begin{center}
\includegraphics[angle=-60,scale=1]{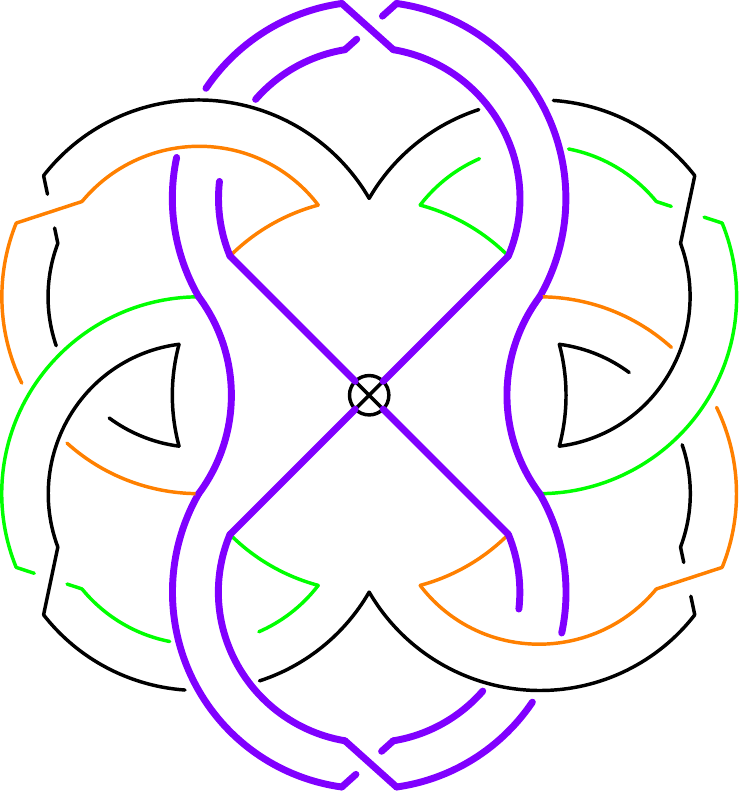}
\caption{An invariant surface of $\Phi_3$ whose dilatation is the golden ratio.}
\label{fig:invariant}
\end{center}
\end{figure}

Hence, the golden ratio which is $\frac{1+\sqrt{5}}{2}$ and which is a root of $x^2-x-1$, is an eigenvalue of $M_{22}$ and hence an eigenvalue of $M_{3}$.
Note that the subsurface generated by $c_3$ and $c_6$ is invariant under the mapping class $\Phi_3$. The invariant subsurface is enclosed by the blue violet curve which is simple in Figure~\ref{fig:invariant}. On that surface, the dilatation is the golden ratio.

\subsection{proof of the Theorem~\ref{thm:main} on $\Sigma_{2k,k}$}
On $\Sigma_{2k,k}$, $r^{k-1}$ sends $(c_{k},c_{2k})$ to $(c_{1},c_{k+1})$. By applying $T_{c_1}$ to $c_i$'s, all $c_i$'s which intersect $c_1=r^{k-1}(c_k)$ wrap around not only $c_i$ but also $c_{1}$. In particular, $c_{k+1}$ wraps around not only $c_{k+1}$ but also $c_{1}$. 
Now, $r$ sends $(c_{1},c_{k+1})$ to $(c_{2k},c_{k})$. 
Hence, $\Phi_k$ is represented by

$$\begin{bmatrix}
0 & I \\
I & 0 \\
\end{bmatrix}
+A$$
where $A=\begin{bmatrix} a_{ij} \end{bmatrix}$ with $a_{ij}=0$ if $i \neq 2k$ and $a_{2kj}=\text{card}(c_k \cap c_j)$. That is

\[\left[\begin{array}{ccc|c|c !{\color{black}\vrule width1pt}ccc|c|c}
& & & & & & & & &  \\
& \mathbf{0} & & \mathbf{0} & \mathbf{0} & & I & &\mathbf{0} & \mathbf{0} \\
& & & & & & & & &  \\
\hline
0 & \cdots & 0 & 0 & 0 &0 & \cdots & 0  & 1 & 0 \\
\hline
0 & \cdots & 0 & 0 & 0 &0 & \cdots & 0 &  0 & 1 \\%
\noalign{
\color{black}
\hrule height 1pt
}%
& & & & & & & & &  \\
& I & & \mathbf{0} & \mathbf{0} & & \mathbf{0} & &\mathbf{0} & \mathbf{0} \\
& & & & & & & & &  \\
\hline
0 & \cdots & 0 & 1 & 0 & 0&  \cdots & 0 & 0 & 0 \\
\hline
 &  \ast  &  & 0 & 1 & &  \ast &  & 1 & 1 \\
\end{array}\right].
\]

By changing the order of basis by switching $\mu_k$ and $\mu_{2k-1}$, $\Phi_k$ can be represented by the following reducible matrix $M_k$,

\[\left[\begin{array}{ccc|c|c|ccc !{\color{blueviolet}\vrule width1pt}c|c}
& & & & & & & & &  \\
& \mathbf{0} & & \mathbf{0} & \mathbf{0} & & I & &\mathbf{0} & \mathbf{0} \\
& & & & & & & & &  \\
\hline
0 & \cdots & 0 & 0 & 1 &0 & \cdots & 0  & 0 & 0 \\
\hline
0 & \cdots & 0 & 1 & 0 &0 & \cdots & 0 &  0 & 0 \\
\hline
& & & & & & & & &  \\
& I & & \mathbf{0} & \mathbf{0} & & \mathbf{0} & &\mathbf{0} & \mathbf{0} \\
& & & & & & & & &  \\%
\noalign{
\color{blueviolet}
\hrule height 1pt
}%
0 & \cdots & 0 & 0 & 0 & 0&  \cdots & 0 & 0 & 1 \\
\hline
 &  \ast  &  & 0 & 1 & &  \ast &  & 1 & 1 \\
\end{array}\right].
\]

$M_k$ can be partitioned into four submatrices. The partitioned matrix can be written as
$$\begin{bmatrix}
M_{11} & \mathbf{0} \\
M_{21}& M_{22} \\
\end{bmatrix}.$$

where $M_{11}$ is the $2(k-1) \times 2(k-1)$ matrix, $\mathbf{0}$ is the $2(k-1) \times 2$ matrix, $M_{21}$ is the $2 \times 2(k-1)$ matrix, and
$M_{22}$ is the $2 \times 2$ matrix.

The characteristic polynomial of $M_k$ is 
\begin{align*}
\text{det}(xI-M_{k} )&=
\text{det}
\begin{bmatrix}
xI-M_{11} & \mathbf{0} \\
\mathbf{0} & I \\
\end{bmatrix}
\text{det}
\begin{bmatrix}
I & \mathbf{0}\\
-M_{21} & I \\
\end{bmatrix}
\text{det}
\begin{bmatrix}
I & \mathbf{0} \\
\mathbf{0} & xI-M_{22}\\
\end{bmatrix} \\
&=\text{det}(xI-M_{11} )\text{det}(xI-M_{22} ) \\
&=(x+1)^{k-1}(x-1)^{k-1}(x^2-x-1).
\end{align*}

Hence, the golden ratio which is $\frac{1+\sqrt{5}}{2}$ and which is a root of $x^2-x-1$, is an eigenvalue of $M_{22}$ and hence an eigenvalue of $M_{k}$.
Note that the subsurface generated by $c_{k}$ and $c_{2k}$ is invariant under the mapping class $\Phi_k$. On that surface, the dilatation is the golden ratio.
\section{acknowledgement}
We thank Erwan Lanneau, Livio Liechti, Julien Marché, Alan Reid, and Darren Long. This work was supported by Basic Science Research Program through the National Research Foundation of Korea (NRF) funded by the Ministry of Education, Science and Technology (No. NRF-2018005847). The second author was supported by 2017 Hongik University Research Fund.

\begin{thebibliography}{Thu88}

\bibitem[CB88]{CassonBleiler88}
Andrew~J. Casson and Steven~A. Bleiler, \emph{Automorphisms of surfaces after
  {N}ielsen and {T}hurston}, London Mathematical Society Student Texts, vol.~9,
  Cambridge University Press, Cambridge, 1988. \MR{964685}

\bibitem[LS18a]{LiechtiStrenner18-1}
Livio Liechti and Bal\'{a}zs Strenner, \emph{The arnoux-yoccoz mapping classes
  via penner's construction}, \url{arXiv:1805.01248}, 2018, Preprint.

\bibitem[LS18b]{LiechtiStrenner18}
\bysame, \emph{Minimal pseudo-{A}nosov stretch factors on nonorientable
  surfaces}, \url{arXiv:1806.00033}, 2018, Preprint.

\bibitem[Pen88]{Penner88}
Robert~C. Penner, \emph{A construction of pseudo-{A}nosov homeomorphisms},
  Trans. Amer. Math. Soc. \textbf{310} (1988), no.~1, 179--197. \MR{930079}

\bibitem[Pen91]{Penner91}
R.~C. Penner, \emph{Bounds on least dilatations}, Proc. Amer. Math. Soc.
  \textbf{113} (1991), no.~2, 443--450. \MR{1068128}

\bibitem[Thu88]{Thurston88}
William~P. Thurston, \emph{On the geometry and dynamics of diffeomorphisms of
  surfaces}, Bull. Amer. Math. Soc. (N.S.) \textbf{19} (1988), no.~2, 417--431.
  \MR{956596}

\end{thebibliography}
\providecommand{\bysame}{\leavevmode\hbox to3em{\hrulefill}\thinspace}
\providecommand{\MR}{\relax\ifhmode\unskip\space\fi MR }
\providecommand{\MRhref}[2]{%
  \href{http://www.ams.org/mathscinet-getitem?mr=#1}{#2}
}
\providecommand{\href}[2]{#2}

\end{document}